\documentclass[11pt,a4paper]{amsart}
\usepackage{graphicx}
\usepackage{amssymb}
\usepackage{amsmath}
\usepackage{amsthm}
\usepackage{amscd}
\usepackage[all,2cell,dvips]{xy}
 \UseAllTwocells \SilentMatrices
\newtheorem{thm}{Theorem}[section]

\newtheorem{lem}[thm]{Lemma}

\newtheorem{prop}[thm]{Proposition}
\theoremstyle{definition}

\theoremstyle{remark}
\newtheorem{rem}[thm]{\bf Remark}
\numberwithin{equation}{section}

\begin{document}
\title[Relative Singularity Categories and Gorenstein-Projective Modules]
{Relative Singularity Categories and Gorenstein-Projective Modules}
\author[X. W. Chen] {Xiao-Wu Chen}
%\thanks{$^*$ The corresponding author}
\thanks{Supported in part by the National Nature Science Foundation of China (Grant
No. 10501041, 10601052).}
%%\subjclass{}%
\thanks{E-mail: xwchen$\symbol{64}$mail.ustc.edu.cn}
%\keywords{Quantum groups, quivers}%
\maketitle
\date{}%
\dedicatory{}%
\commby{}%

\begin{center}
Department of Mathematics \\
University of Science and Technology of China \\
Hefei 230026, P. R. China
\end{center}

\begin{abstract}

We introduce the notion of relative singularity category with
respect to any self-orthogonal subcategory $\omega$ of an abelian
category. We introduce the Frobenius category of
$\omega$-Cohen-Macaulay objects, and under some reasonable
conditions, we show that the stable category of
$\omega$-Cohen-Macaulay objects is triangle-equivalent to the
relative singularity category. As applications, we relate the stable
category of (unnecessarily finitely-generated) Gorenstein-projective
modules with  singularity categories of rings. We prove that for a
Gorenstein ring, the stable category of Gorenstein-projective
modules is compactly generated and its compact objects coincide with
finitely-generated Gorenstein-projective modules up to direct
summands.
\end{abstract}

\section{Introduction}

\subsection{}

Throughout, $\mathcal{A}$ is an abelian category, $\omega\subseteq
\mathcal{A}$ its full additive subcategory. Denote by
$C^b(\mathcal{A})$, $K^b(\mathcal{A})$ and $D^b(\mathcal{A})$ the
category of bounded complexes, the bounded homotopy category and the
bounded derived category of $\mathcal{A}$, respectively,  both of
whose shift functors will be denoted by [1]. Recall that for any $X,
Y \in \mathcal{A}$, the $n$-th extension group ${\rm
Ext}_\mathcal{A}^n(X, Y)$ is defined to be ${\rm
Hom}_{D^b(\mathcal{A})}(X, Y[n])$, $n\geq 0$ (see \cite{Har}, p.62).
The subcategory $\omega$ is said to be \emph{self-orthogonal} if for
any $X, Y \in \omega$, $n\geq 1$, ${\rm Ext}_\mathcal{A}^n(X, Y)=0$.
Consider the following composite of functors
\begin{align*}
K^b(\omega) \longrightarrow K^b(\mathcal{A}) \longrightarrow
D^b(\mathcal{A}),
\end{align*}
where the first is the inclusion functor, and the second the
quotient functor. By \cite{Ha1}, Chapter II, Lemma 3.4 (or Chapter
III, Lemma 2.1), the composite functor is fully-faithful if and only
if $\omega$ is self-orthogonal.

Let $\omega\subseteq \mathcal{A}$ be a self-orthogonal additive
subcategory. By the argument above, we may view $K^b(\omega)$ as a
triangulated subcategory of $D^b(\mathcal{A})$. Define the
\emph{relative singularity category} $D_\omega(\mathcal{A})$ of
$\mathcal{A}$ with respect to $\omega$ to be the Verdier quotient
category
\begin{align*}
D_\omega(\mathcal{A}): = D^b(\mathcal{A})/{K^b(\omega)}.
\end{align*}

The motivation of introducing relative singularity category is
twofold: (1) A special case of relative singularity category is of
particular interest: let $\mathcal{P}$ denote the subcategory
consisting of projective objects, which is clearly self-orthogonal,
then the relative singularity category with respect to $\mathcal{P}$
is called the \emph{singularity category} of $\mathcal{A}$. Denote
it by $D_{\rm sg}(\mathcal{A})$ (compare \cite{O1}). This
terminology is justified by the fact that: the singularity category
$D_{\rm sg}(\mathcal{A})$ vanishes if and only if the category
$\mathcal{A}$ has enough projectives and every object is of finite
projective dimension. (2) Singularity categories may be described as
relative singularity categories via tilting subcategories, and this
viewpoint allows us to describe singularity category by various
tilting subcategories. Precisely, let $\mathcal{A}$ have enough
projectives, a \emph{tilting subcategory} $\mathcal{T}$ is a
self-orthogonal subcategory such that
$K^b(\mathcal{T})=K^b(\mathcal{P})$ inside $D^b(\mathcal{A})$. Then
we have $D_{\rm sg}(\mathcal{A})=D_{\mathcal{T}}(\mathcal{A})$ for
any tilting subcategory $\mathcal{T}$.

\subsection{} The paper is organized as follows: In section 2, we
study the relative singularity category $D_{\omega}(\mathcal{A})$
and various related subcategories of the abelian category
$\mathcal{A}$, in particular, the category of
$\omega$-Cohen-Macaulay objects. As the main theorem, we prove that
there is a full exact embedding of the stable category of
$\omega$-Cohen-Macaulay objects into the relative singularity
category, and further under some reasonable conditions, the
embedding is an equivalence. In section 3, we apply the result to
the module category of rings, and we rediscover the result of
Buchweitz-Happel which says that for a Gorenstein ring, the
singularity category is triangle-equivalent to the stable category
of finitely-generated Gorenstein-projective modules, and we also
find a similar result holds in the unnecessarily finitely-generated
case. We relate the stable category of  $T$-Cohen-Macaulay objects
to the stable category of      finitely-generated
Gorenstein-projective modules over the endomorphism ring ${\rm
End}_\mathcal{A}(T)$, where $T$ is any self-orthogonal object in
$\mathcal{A}$. In section 4, we show that for a Gorenstein ring, the
stable category of Gorenstein-projective modules is compactly
generated and its subcategory of compact objects is the stable
category of finitely-generated Gorenstein-projective modules up to
direct summands.

For triangulated categories, we refer to \cite{Ha1, Har, V1}. We
abuse the notions of triangle-functors and exact functors between
triangulated categories. For Gorenstein rings and
Gorenstein-projective modules, we refer to \cite{EJ, Buc, Ha2, BR}.

\section{Relative singularity category and $\omega$-Cohen-Macaulay objects}

\subsection{} In this subsection, we will introduce some
subcategories of the abelian category $\mathcal{A}$ (compare
\cite{AR,CZ}). At this moment,  $\omega\subseteq \mathcal{A}$ is an
arbitrary additive subcategory. Consider the following full
subcategories:
\begin{align*}
\widehat{\omega}:= \{ X \in \mathcal{A} \; &| \; \mbox{ there exists
an exact sequence }\\
&\; 0 \rightarrow T^{-n} \rightarrow T^{1-n} \rightarrow \cdots
\rightarrow T^0 \rightarrow X \rightarrow 0,
\mbox{ each } \; T^{-i}\in \omega,\;  n\geq 0 \};\\
\omega^\perp := \{X \in \mathcal{A}\; &|\; \; {\rm
Ext}_\mathcal{A}^i(T, X)=0, \;
\mbox{ for all } \; T \in \omega,\; i\geq 1\};\\
_\omega \mathcal{X} : =\{X \in \mathcal{A}\; & |\; \mbox{ there exists an exact sequence } \\
&\dots  \rightarrow T^{-n} \stackrel{d^{-n}}{\rightarrow }T^{1-n}
\rightarrow \cdots \rightarrow T^0 \stackrel{d^0}{\rightarrow} X
\rightarrow 0, \mbox{ each } \; T^{-i}\in \omega,  \; {\rm
Ker}d^i\in \omega^\perp \}.
\end{align*}
If $\omega$ is self-orthogonal, using the dimension-shift technique
in homological algebra, we infer that $\widehat{\omega}\subseteq
\omega^\perp$ and $_\omega \mathcal{X} \subseteq \omega^\perp$, and
thus we get $\widehat{\omega}\subseteq {_\omega \mathcal{X}}$.
Consequently, if $\omega$ is self-orthogonal, we obtain that
\begin{align*}
\omega \subseteq \widehat{\omega} \subseteq {_\omega \mathcal{X}}
\subseteq {\omega^\perp}.
\end{align*}

Dually, we have the following three full subcategories:
\begin{align*}
\stackrel{ \vee}{\omega}:= \{ X \in \mathcal{A} \; &| \; \mbox{
there exists an exact sequence }\\
& \; 0  \rightarrow X \rightarrow T^0 \rightarrow \cdots \rightarrow
T^{n-1} \rightarrow T^n \rightarrow 0,
\mbox{ each   } \; T^{i}\in \omega, \; n\geq 0 \};\\
{^\perp \omega} := \{X \in \mathcal{A}\; &|\; \; {\rm
Ext}_\mathcal{A}^i(X, T)=0, \;
\mbox{ for all } \; T \in \omega, \; i\geq 1\};\\
\mathcal{X}_\omega :=\{X \in \mathcal{A}\; & |\; \mbox{ there exists an exact sequence } \\
 0 \rightarrow & X \stackrel{d^{-1}}{\rightarrow} T^0
\stackrel{d^0}{\rightarrow }T^1 \rightarrow \cdots \rightarrow
T^{n-1} \stackrel{d^{n-1}}{\rightarrow} T^n \rightarrow \cdots,
\mbox{ each   } \; T^{i}\in \omega,\;  {\rm Coker}d^i\in {^\perp
\omega} \}.
\end{align*}
Similarly, if $\omega$ is self-orthogonal, we have $\omega \subseteq
\stackrel{\vee}{\omega} \; \subseteq {\mathcal{X}_\omega} \subseteq
{^\perp \omega} $.

Let $\omega$ be a self-orthogonal subcategory. We define the
\emph{category of $\omega$-Cohen-Macaulay objects} to be the full
subcategory $\alpha(\omega):= {\mathcal{X}_\omega} \cap {_\omega
\mathcal{X}}$. By \cite{AR}, Proposition 5.1, the full subcategories
$_\omega \mathcal{X}$ and $\mathcal{X}_\omega$ are closed under
extensions, and therefore so is $\alpha(\omega)$. Hence,
$\alpha(\omega)$ becomes an exact category whose conflations are
just short exact sequences with terms in $\alpha(\omega)$ (for
terminology, see \cite{Ke3}). Observe that objects in $\omega$ are
(relative) projective and injective, and then it is not hard to see
that $\alpha(\omega)$ is a Frobenius category, whose
projective-injective objects are precisely contained in the additive
closure ${\rm add}\; \omega$ of $\omega$. Consider the stable
category $\underline{\alpha(\omega)}$ of $\alpha(\omega)$ modulo
$\omega$ (or equivalently modulo ${\rm add}\; \omega$). Then by
\cite{Ha1}, $\underline{\alpha(\omega)}$ is a triangulated category.

For each $X\in \alpha(\omega)$, from the definition (and the
dimension-shift technique if needed), we have an exact sequence in
$K(\omega)$
$$T^\bullet= \quad \rightarrow \cdots \rightarrow T^{-n} \rightarrow T^{-n+1} \rightarrow \cdots \rightarrow
T^{-1} \rightarrow T^0 \rightarrow T^1 \rightarrow \cdots
\rightarrow T^n\rightarrow T^{n+1} \rightarrow \cdots$$ such that
each of its cocycles $Z^i(T^\bullet)$ lies in ${^\perp \omega}\cap
\omega^\perp$, and $X=Z^0(T^\bullet)$. Such a complex $T^\bullet$
will be called an \emph{$\omega$-complete resolution} for $X$. It is
worthy observing that an exact complex $T^\bullet\in K(\omega)$ is
an $\omega$-complete resolution if and only if for each $T\in
\omega$, the Hom complexes ${\rm Hom}(T, T^\bullet)$ and ${\rm
Hom}(T^\bullet, T)$ are exact. One may compare \cite{Bel},
Definition 5.5.

\subsection{}

Consider the following composite of natural functors
\begin{align*}
F: \alpha(\omega)  \longrightarrow \mathcal{A}
\stackrel{i_\mathcal{A}}\longrightarrow D^b(\mathcal{A})
\stackrel{Q_\omega}{\longrightarrow} D_\omega(\mathcal{A}),
\end{align*}
where the first functor is the inclusion, the second is the full
embedding which sends objects in $\mathcal{A}$ to the stalk
complexes concentrated at degree $0$, and the last is the quotient
functor $Q_\omega: D^b(\mathcal{A})\longrightarrow
D_{\omega}(\mathcal{A})$. Note that $F(\omega)=0$, and thus $F$
induces a unique functor $\underline{F}$ from
$\underline{\alpha(\omega)}$ to $D_\omega(\mathcal{A})$.

Our main result is

\begin{thm}
Let $\omega \subseteq \mathcal{A}$ be a self-orthogonal additive
subcategory. Then the natural functor $\underline{F}:
\underline{\alpha(\omega)} \longrightarrow D_\omega(\mathcal{A})$ is
a fully-faithful triangle-functor. \par

Assume further that $\widehat{\mathcal{X}_\omega}= \mathcal{A}
=\stackrel{\vee}{_\omega \mathcal{X}}$. Then $\underline{F}$ is an
equivalence, thus a triangle-equivalence.
\end{thm}

Note that the subcategories $\widehat{\mathcal{X}_\omega}$ and
$\stackrel{\vee}{_\omega \mathcal{X}}$ are defined as in {\bf 2.1},
by replacing $\omega$ by $\mathcal{X}_\omega$ and
$_\omega\mathcal{X}$, respectively.

\subsection{} We will divide the proof of Theorem 2.1 into proving several
propositions. Note that we will always view $\mathcal{A}$ as the
full subcategory of $D^b(\mathcal{A})$ consisting of stalk complexes
concentrated at degree zero.

We need some notation. A complex $X^\bullet \in C^b(\omega)$ is said
to \emph{negative} if $X^n=0$ for all $n\geq 0$. Denote by
$D^{<0}(\omega)$ to be full subcategory of $K^b(\omega)$ whose
objects are isomorphic t o some negative complexes in $C^b(\omega)$.
Similarly, we have the subcategory $D^{>0}(\omega)$.

\begin{lem}
(1).\quad For $M\in {^\perp\omega}$ and $X^\bullet\in D^{<
0}(\omega)$, we have ${\rm Hom}_{D^b(\mathcal{A})}(M, X^\bullet)=0$.\\
(2).\quad For $N\in \omega^\perp$ and $Y^\bullet\in D^{>
0}(\omega)$, we have ${\rm Hom}_{D^b(\mathcal{A})}(Y^\bullet, N)=0$.
\end{lem}

\noindent {\bf Proof.}\quad We only show (1). Consider
$\mathcal{L}:=\{Z^\bullet\in D^b(\mathcal{A})\; |\; {\rm
Hom}_{D^b(\mathcal{A})}(M, Z^\bullet)=0\}$. By the self-orthogonal
property of $\omega$, we have $\omega[i]\in \mathcal{L}$ for all
$i>0$. Observe that $\mathcal{L}$ is closed under extensions, and
complexes in $D^{<0}(\omega)$ are obtained by iterated extensions
from objects in $\bigcup_{i>0}\omega[i]$, thus we infer that
$D^{<0}(\omega)\subseteq \mathcal{L}$. \hfill $\blacksquare$

In what follows, morphisms in $D^b(\mathcal{A})$ will be denoted by
arrows, and those whose cones lie in $K^b(\omega)$ will be denoted
by doubled arrows; morphisms in $D_\omega(\mathcal{A})$ will be
denoted by right fractions (for the definition, see \cite{V1}).

Let $M, N \in \mathcal{A}$. We consider the natural map
$$\theta_{M, N}: \; {\rm Hom}_{\mathcal{A}}(M, N) \longrightarrow
{\rm Hom}_{D_\omega(\mathcal{A})}(Q_\omega(M), Q_\omega(N)), \quad
f\longmapsto f/{{\rm Id}_M}.
$$
Set $\omega(M, N)=\{f\in {\rm Hom}_\mathcal{A}(M, N)\; |\; f \mbox{
factors through objects in }\omega\}$. Then $\theta_{M, N}$ vanishes
on $\omega(M, N)$ because $Q_\omega(\omega)=0$.

The following observation is crucial in our proof, compare
\cite{CZ}, Lemma 2.1 and \cite{O1}, Proposition 1.21.

\begin{lem}
In the following two cases: (1) $M \in \mathcal{X}_\omega$ and $N\in
\omega^\perp$; (2) $M \in {^\perp \omega}$ and $N\in {_\omega
\mathcal{X}}$, the morphism $\theta_{M,N}$ induces an isomorphism
$${\rm Hom}_\mathcal{A}(M, N)/{\omega(M, N)}\simeq {\rm Hom}_{D_{\omega}(\mathcal{A})}(M, N).$$
\end{lem}

\noindent {\bf Proof.}\quad  We only show case (1). First, we show
that $\theta_{M, N}$ is surjective. For this, consider any morphism
$a/s: M \stackrel{\; s}\Longleftarrow Z^\bullet
\stackrel{a}\longrightarrow N$ in $D_{\omega}(\mathcal{A})$, where
$Z^\bullet$ is a complex, both $a$ and $s$ are morphisms in
$D^b(\mathcal{A})$, and the cone of $s$, $C^\bullet={\rm Con}(s)$,
lies in $K^b(\omega)$. Hence we have a distinguished triangle in
$D^b(\mathcal{A})$
\begin{align}\label{eqn2.1}
Z^\bullet \stackrel{s}{\Longrightarrow } M \longrightarrow C^\bullet
\longrightarrow Z^\bullet[1].
\end{align}

Since $M\in \mathcal{X}_\omega$, we have a long exact sequence
\begin{align*}
0 \longrightarrow M \stackrel{\varepsilon}\longrightarrow T^0
\stackrel{d^0}\longrightarrow
 T^1 \longrightarrow \cdots \longrightarrow T^n \stackrel{d^n}\longrightarrow T^{n+1} \longrightarrow \cdots
\end{align*}
where each $T^i\in \omega$ and ${\rm Ker}d^i \in {^\perp \omega}$.
Hence in $D^b(\mathcal{A})$, $M$ is isomorphic to the following
complex
$$T^\bullet:= \;\;  0 \longrightarrow T^0 \stackrel{d^0}\longrightarrow T^1 \longrightarrow \cdots
\longrightarrow T^n \stackrel{d^n}\longrightarrow T^{n+1}
\longrightarrow  \cdots, $$
 and furthermore, $M$ is isomorphic to
the good truncation $\tau^{\leq l}T^\bullet$ for any $l\geq 0$. Note
the following natural triangle in $K^b(\mathcal{A})$
\begin{align}\label{eqn2.2}
(\sigma^{< l}T^\bullet) [-1] \longrightarrow {\rm Ker}d^l[-l]
\stackrel{s''}{\Longrightarrow} \tau^{\leq l}T^\bullet
\longrightarrow \sigma^{< l}T^\bullet,
\end{align}
where $\sigma^{<l}T^\bullet$ is the brutal truncation. Take $s'$ to
be the following composite in $D^b(\mathcal{A})$
$${\rm Ker}d^l[-l] \stackrel{s''}\Longrightarrow \tau^{\leq l} T^\bullet \longrightarrow T^\bullet
\stackrel{\varepsilon}\Longleftarrow  M.$$
Thus from the triangle
(\ref{eqn2.2}), we get a triangle in $D^b(\mathcal{A})$
\begin{align}\label{eqn2.3}
(\sigma^{< l} T^\bullet)[-1] \longrightarrow {\rm Ker}d^l[-l]
\stackrel{s'}\Longrightarrow M \stackrel{\varepsilon}\longrightarrow
\sigma^{< l}T^\bullet.
\end{align}

Since $C^\bullet\in K^b(\omega)$, we may assume that
$$C^\bullet = \quad \cdots \longrightarrow 0 \longrightarrow W^{-t'} \longrightarrow
\cdots \longrightarrow W^{t-1}\longrightarrow W^t \longrightarrow 0
\longrightarrow \cdots,$$ where $W^i\in \omega$, $t, t'\geq 0$. Set
$l_0=t+1$, $E={\rm Ker}d^{l_0}$. Note that $E\in {^\perp \omega}$
and $C^\bullet [l_0]\in D^{< 0}(\omega)$, by Lemma 2.2(1), we get
$${\rm Hom}_{D^b(\mathcal{A})}(E[-l_0], C^\bullet)= {\rm
Hom}_{D^b(\mathcal{A})}(E, C^\bullet[l_0])=0.$$
Hence, the morphism
$E[-l_0]\stackrel{s'} \Longrightarrow M \longrightarrow C^\bullet$
is $0$. By the triangle (\ref{eqn2.1}), we infer that there exists
$h: E[-l_0]\longrightarrow Z^\bullet$ such that $s'=s\circ h$, and
thus $a/s=(a\circ h)/{s'}$.

Note that $N \in \omega^\perp$ and $(\sigma^{< l_0}T^\bullet)[-1]\in
D^{>0}(\omega)$, by Lemma 2.2(2), we have
$${\rm
Hom}_{D^b(\mathcal{A})}((\sigma^{< l_0}T^\bullet)[-1], N)=0.$$
Applying the cohomological functor ${\rm Hom}_{D^b(\mathcal{A})}(- ,
N)$ to the triangle (\ref{eqn2.3}), we obtain the following exact
sequence (here, take $l=l_0$)
\begin{align*}
{\rm Hom}_{D^b(\mathcal{A})}(M, N) \stackrel{{\rm
Hom}_{D^b(\mathcal{A})}(s', N)} \longrightarrow {\rm
Hom}_{D^b(\mathcal{A})}(E[-l_0], N) \longrightarrow {\rm
Hom}_{D^b(\mathcal{A})}((\sigma^{< l_0}T^\bullet)[-1], N).
\end{align*}
Thus, there exists $f: M \longrightarrow N$ such that $f\circ
s'=a\circ h$. Hence, we have
$$a/s= (a\circ h)/s'= (f\circ s')/s'=\theta_{M, N}(f),$$
proving that $\theta_{M, N}$ is surjective.

Next, we will show ${\rm Ker}\theta_{M, N}=\omega(M, N)$, then we
are done. It is already known that $\omega(M, N)\subseteq {\rm
Ker}\theta_{M, N}$. Conversely, consider $f: M \longrightarrow N$
such that $\theta_{M, N}(f)=0$. Hence there exists $s: Z^\bullet
\Longrightarrow M$ such that $f\circ s=0$, where $s$ is a morphism
in $D^b(\mathcal{A})$ whose cone $C^\bullet \in K^b(\omega)$. Using
the notation above, we obtain that $s'=s\circ h$. Thus $f\circ
s'=0$. By the triangle (\ref{eqn2.3}), we infer that there exists
$f': \sigma^{< l_0}T^\bullet \longrightarrow N$  such that $f'\circ
\varepsilon=f$.

Consider the following natural triangle
\begin{align}\label{eqn2.4}
T^0[-1] \longrightarrow \sigma^{>0}(\sigma^{<l_0}T^\bullet)
\Longrightarrow \sigma^{< l_0}T^\bullet
\stackrel{\pi}\longrightarrow T^0.
\end{align}
Since $N \in \omega^\perp$ and
$\sigma^{>0}(\sigma^{<l_0}T^\bullet)\in D^{>0}(\omega)$, by Lemma
2.2(2), we have
$${\rm
Hom}_{D^b(\mathcal{A})}(\sigma^{>0}(\sigma^{<l_0}T^\bullet), N)=0.$$
Thus the composite morphism
$\sigma^{>0}(\sigma^{<l_0}T^\bullet)\Longrightarrow
\sigma^{<l_0}T^\bullet \stackrel{f'} \longrightarrow N$ is $0$, and
furthermore, by the triangle (\ref{eqn2.4}), we infer that there
exists $g: T^0\longrightarrow N$ such that $g\circ \pi=f'$. So we
get $f=g\circ (\pi\circ \varepsilon)$, which proves that $f$ factors
through $\omega$ inside $D^b(\mathcal{A})$. Note again that
$i_\mathcal{A}: \mathcal{A} \longrightarrow D^b(\mathcal{A})$ is
fully-faithful, and we can obtain that $f$ factors through $\omega$
in $\mathcal{A}$, i.e., $f\in \omega(M, N)$. This finishes the
proof. \hfill $\blacksquare$

Recall the notion of $\partial$-functor, compare \cite{Ke}, section
1. Let $(\mathfrak{a}, \; \mathcal{E})$ be an exact category,
$\mathcal{C}$ a triangulated category. An additive functor $F:
\mathfrak{a} \longrightarrow \mathcal{C}$ is said to be a
$\partial$-functor, if for each conflation $(i, \; d):\; X
\stackrel{i}{\longrightarrow} Y \stackrel{d}\longrightarrow Z \; \in
\mathcal{E}$, there exists a morphism $w_{(i, d)}: F(Z)
\longrightarrow F(X)[1]$ such that the triangle is distinguished
$$F(X) \stackrel{F(i)}\longrightarrow F(Y) \stackrel{F(d)}\longrightarrow F(Z)
 \stackrel{w_{(i, d)}}\longrightarrow F(X)[1],$$
moreover, the morphisms $w$ are natural in the sense that given a
morphism between two conflations
\[
\xymatrix{
X \ar[r]^{i} \ar[d]^-{f} & Y \ar[d]^-g \ar[r]^-{d}  & Z \ar[d]^-h\\
X' \ar[r]^-{i'} & Y' \ar[r]^-{d'} &Z', }\]
 then we have a morphism of triangles
\[
\xymatrix{
F(X) \ar[r]^{F(i)} \ar[d]^-{F(f)} & F(Y) \ar[d]^-{F(g)} \ar[r]^-{F(d)}  & F(Z)
\ar[d]^-{F(h)} \ar[r]^-{w_{(i,\; d)}} & F(X)[1] \ar[d]^-{F(f)[1]}\\
F(X') \ar[r]^-{F(i')} & F(Y') \ar[r]^-{F(d')} &F(Z')
\ar[r]^-{w_{(i', \; d')}} & F(X')[1]. }\]

We will need the following fact, which is direct from definition.

\begin{lem}
Let $F: \mathfrak{a}\longrightarrow \mathcal{C}$ be a
$\partial$-functor. Assume $j: \mathfrak{b}\longrightarrow
\mathfrak{a}$ is an exact functor between two exact categories,
$\pi: \mathcal{C} \longrightarrow \mathcal{D}$ a triangle-functor
between two triangulated categories. Then the composite functor $\pi
Fj:\mathfrak{b}\longrightarrow \mathcal{D}$ is a $\partial$-functor.
\end{lem}

Next fact is very useful, and well-known, compare \cite{Ha1}, p.23.

\begin{lem}
Let $(\mathfrak{a},\; \mathcal{E})$ be a Frobenius category,
$\underline{\mathfrak{a}}$ its stable category modulo projectives.
Assume $F: \mathfrak{a}\longrightarrow \mathcal{C}$ is a
$\partial$-functor, which vanishes on projective objects. The
induced functor $\underline{F}: \underline{\mathfrak{a}}
\longrightarrow \mathcal{C}$ is a triangle-functor.
\end{lem}

\noindent {\bf Proof.}\quad Since $F$ vanishes on projective
objects, then the functor $\underline{F}$ is defined. Recall that
the translation functor $S$ on $\underline{\mathfrak{a}}$ is defined
such that for each $X$, we have a fixed conflation
$X\stackrel{i_X}\longrightarrow I(X) \stackrel{d_X} \longrightarrow
S(X)$, where $I(X)$ is injective (for details, see \cite{Ha1}). By
assumption, we have the distinguished triangle in $\mathcal{C}$
$$F(X) \stackrel{F(i_X)}\longrightarrow F(I(X)) \stackrel{F(d_X)}\longrightarrow F(S(X))
\stackrel{w_{(i_X, d_X)}}\longrightarrow F(X)[1].$$
 Since
$F(I(X))\simeq 0$, we infer that $w_{(i_X, d_X)}$ is an isomorphism.
Set $\eta_X:=w_{(i_X, d_X)}$. In fact, by the naturalness of $w$, we
can obtain that $\eta_X$ is natural in $X$, in other words, $\eta:
FS\longrightarrow [1]F$ is a natural isomorphism. Recall that all
the distinguished triangles in $\underline{\mathfrak{a}}$ arise from
conflations in $\mathfrak{a}$ (\cite{CZ}, Lemma 2.1), then one may
show that $(\underline{F}, \eta)$ is a triangle-functor easily. We
omit the details. \hfill $\blacksquare$

\vskip 10pt

\noindent {\bf Proof of Theorem 2.1:}\quad By Lemma 2.3, we know
that $\underline{F}$ is fully-faithful. It is classical that
$i_\mathcal{A}: \mathcal{A} \longrightarrow D^b(\mathcal{A})$ is a
$\partial$-functor (by \cite{Har}, p.63, Remark). Then by Lemma 2.4,
we know that the composite functor $F$ is also a $\partial$-functor.
Now by Lemma 2.5, we deduce that $\underline{F}$ is a
triangle-functor.

Now assume that $\widehat{\mathcal{X}_\omega} = \mathcal{A}
=\stackrel{\vee}{_\omega \mathcal{X}}$. It suffices to show that
$\underline{F}$ is dense, that is, the image ${\rm
Im}\underline{F}=D_\omega(\mathcal{A})$. By above, we know that
${\rm Im}\underline{F}$ is a triangulated subcategory, and it is
direct to see that $D_{\omega}(\mathcal{A})$ is generated by the
image $Q_\omega(\mathcal{A})$ of $\mathcal{A}$ in the sense of
\cite{Ha1}, p.71. Hence it is enough to show that
$Q_\omega(\mathcal{A})$ lies in ${\rm Im}\underline{F}$.

Assume $X\in \mathcal{A}$. Since $\omega$ cogenerates
$\mathcal{X}_\omega$ and $X\in
\widehat{\mathcal{X}_\omega}=\mathcal{A}$, by Auslander-Buchweitz
decomposition theorem (\cite{AB}, Theorem 1.1), we have an exact
sequence
$$0 \longrightarrow Y \longrightarrow X' \longrightarrow X\longrightarrow 0,$$
where $Y\in \widehat{\omega}$, and $X'\in \mathcal{X}_\omega$. Since
$Y\in \widehat{\omega}$, then inside $D^b(\mathcal{A})$ we have
$Y\in K^b(\omega)$. Consequently, $Q_\omega(Y)\simeq 0$. Note that
the above exact sequence induces a distinguished triangle in
$D^b(\mathcal{A})$ (\cite{Har}, p.63), and thus we have the induced
distinguished triangle in $D_\omega(\mathcal{A})$
$$Q_\omega(Y) \longrightarrow Q_\omega(X') \longrightarrow Q_\omega(X) \longrightarrow Q_\omega(Y)[1].$$
Now since $Q_\omega(Y)\simeq 0$, we deduce that $Q_\omega(X')\simeq
Q_\omega(X)$. On the other hand, $\omega$ generates $_\omega
\mathcal{X}$ and $X'\in \; \stackrel{\vee}{_\omega
\mathcal{X}}=\mathcal{A}$, by the dual of Auslander-Buchweitz
decomposition theorem, we have an exact sequence
$$0 \longrightarrow X' \longrightarrow X'' \longrightarrow Y' \longrightarrow 0,$$
where $Y'\in \stackrel{\vee}{\omega}$, and $X''\in
{_\omega\mathcal{X}}$. By the same argument as above, we deduce that
$Q_\omega(X')\simeq Q_\omega(X'')$, and consequently,
$Q_\omega(X)\simeq Q_\omega(X'')$. As we noted in {\bf 2.1} that
$\stackrel{\vee}\omega\subseteq { \mathcal{X}_\omega}$, and in the
exact sequence above, both $Y$ and $X'$ lie in $\mathcal{X}_\omega$,
and by Proposition 5.1 in \cite{AR}, $\mathcal{X}_\omega$ is closed
under extensions, we infer that $X''\in \mathcal{X}_\omega$, and
thus $X''\in \alpha(\omega)$. Hence $Q_\omega(X'')=F(X'')$, and we
see that $Q_\omega(X)$ lies in the image of $\underline{F}$. This
completes the proof. \hfill $\blacksquare$

\section{Gorenstein-projective modules and singularity categories}

\subsection{} Let $R$ be a ring with unit. Denote by $R\mbox{-Mod}$
 the category of left $R$-modules, and $R\mbox{-Proj}$ its full
 subcategory of projective modules. A complex $P^\bullet=(P^n, d^n)$
 in $C(R\mbox{-Proj})$ is said to be \emph{totally-acyclic} (\cite{Kr}, section 7),
  if for each projective module $Q$,
 the Hom complexes ${\rm Hom}_R(Q, P^\bullet)$ and  ${\rm Hom}_R(P^\bullet,
 Q)$ are exact. Hence a complex $P^\bullet$ is totally-acyclic if
 and only if it is acyclic (= exact) and for each $n$, the cocycle ${\rm
 Ker}d^n$ lies in $^\perp{R\mbox{-Proj}}$. A module $M$ is said to
 be \emph{Gorenstein-projective}, if there exists a totally-acyclic complex
 $P^\bullet$ such that its zeroth cocycle is $M$. In this case,
 $P^\bullet$ is said to be a \emph{complete resolution} of $M$.
 Denote by $R\mbox{-GProj}$ the full subcategory consisting of
 Gorenstein-projective modules.

Observe that  a module $M$ is Gorenstein-projective if and only if
there exists an exact sequence $0\longrightarrow M
\stackrel{\varepsilon}\longrightarrow P^0
\stackrel{d^0}\longrightarrow P^1 \stackrel{d^1}\longrightarrow P^2
\longrightarrow \cdots$ such that each cocycle ${\rm Ker}d^i\in
{^\perp R\mbox{-Proj}}$. Set $\mathcal{A}=R\mbox{-Mod}$,
$\omega=R\mbox{-Proj}$. Thus $_\omega \mathcal{X}=\mathcal{A}$ and
$\alpha(\omega)=\mathcal{X}_\omega$. By the above observation, we
have $\alpha(\omega)=R\mbox{-GProj}$. In this case, the relative
singularity category is the (big) singularity category of $R$
(compare \cite{O1})
$$D'_{\rm sg}(R)=D^b(R\mbox{-Mod})/K^b(R\mbox{-Proj}).$$
Note that $D'_{\rm sg}(R)$ vanishes if and only if every module has
finite projective dimension, and then it is equivalent to that the
ring $R$ has finite left global dimension.

The following result can be read from the general theory developed
in section 2.

\begin{prop}
(1)\quad The category $R\mbox{\rm -GProj}$ is a Frobenius category
with projective-injective objects exactly contained in $R\mbox{\rm -Proj}$.\\
(2)\quad The natural functor $\underline{F}: R\underline{\mbox{\rm
-GProj}} \longrightarrow D'_{\rm sg}(R)$ is fully-faithful and
exact.
\end{prop}

A sufficient condition making $\underline{F}$ an equivalence is that
the ring $R$ is Gorenstien. This was first observed by Buchweitz
\cite{Buc}. Recall that a ring $R$ is said to be \emph{Gorenstein},
if $R$ is two-sided noetherian, and the regular module $R$ has
finite injective dimension both as a left and right module.

We need the following fact, which is known to experts.

\begin{lem}
Let $R$ be a Gorenstein ring. Then we have $R\mbox{\rm
-GProj}={^\perp R\mbox{\rm -Proj}}.$
\end{lem}

\noindent {\bf Proof.}\quad Note that $R\mbox{-GProj}\subseteq
{^\perp R\mbox{\rm -Proj}}$. For the converse, denote $\mathcal{L}$
the full subcategory of $R\mbox{-Mod}$, consisting of modules of
finite injective dimension. By \cite{EJ}, Lemma 10.2.13,
$\mathcal{L}$ is preenveloping (= covariantly-finite), i.e., for any
module $M$, there exists a morphism $g_M: \; M \longrightarrow C_M$
such that $C_M\in \mathcal{L}$ and any morphism from $M$ to a module
in $\mathcal{L}$ factors through $g_M$ (such a morphism $g_M$ is
called an $\mathcal{L}$-preenvelop (= right
$\mathcal{L}$-approximation) ). We note that the morphism $g_M$ is
mono, by noting that the injective hull of $M$ factors through
$g_M$.

Now assume $M\in {^\perp R\mbox{-Proj}}$. Take an exact sequence
\begin{align} \label{eqn3.1}
0 \longrightarrow K \longrightarrow P^0
\stackrel{\theta}\longrightarrow C_M\longrightarrow 0,
\end{align}
where $P^0$ is projective. Since $C_M$ has finite injective
dimension, by \cite{EJ}, Proposition 9.1.7, it also has finite
projective dimension. Thus we infer that $K$ has finite projective
dimension. Note that $M\in {^\perp R\mbox{-Proj}}$, and by the
dimension-shift argument, we have ${\rm Ext}_R^1(M, K)=0$. Applying
the functor ${\rm Hom}_R(M, -)$ to (\ref{eqn3.1}), we obtain a long
exact sequence, and from which, we read a surjective map ${\rm
Hom}_R(M, \theta):\; {\rm Hom}_R(M, P)\longrightarrow {\rm Hom}_R(M,
C_M)$. In particular, the morphism $g_M$ factor through $\theta$,
and thus we get a morphism
 $h: M \longrightarrow P^0$ such that $g_M=\theta \circ h$. Since
 $g_M$ is an $\mathcal{L}$-preenvelop, and $g_M$ factors through $h$
 (note $P^0\in \mathcal{L}$), and we deduce that $h$ is also an
 $\mathcal{L}$-preenvelop. In particular, $h$ is mono. Consider the
 exact sequence
\begin{align} \label{eqn3.2}
0 \longrightarrow M \stackrel{h}\longrightarrow P^0 \longrightarrow
M'\longrightarrow 0.
\end{align}
For any projective module $Q$, applying the functor ${\rm Hom}_R(-,
Q)$, and we obtain a long exact sequence, from which we read that
${\rm Ext}^i_R(M, Q)=0$ for $i\geq 1$ (for $i=1$, we need the fact
that  $h$ is an \emph{$\mathcal{L}$-preenvelop}). Thus $M'\in
{^\perp R\mbox{-Proj}}$. Applying the same argument to $M'$, we may
get an exact sequence $0 \longrightarrow M' \longrightarrow
P^1\longrightarrow M''\longrightarrow 0$ with $P$ projective and
$M''\in {^\perp R\mbox{-Proj}}$. Continue this process, and we
obtain a long exact sequence $0\longrightarrow M \longrightarrow P^0
\longrightarrow P^1\longrightarrow P^2 \longrightarrow \cdots$ with
cocycles in ${^\perp R\mbox{-Proj}}$, that is, $M\in
R\mbox{-GProj}$. Thus we are done. \hfill $\blacksquare$

Now we have the following variant of Buchweitz-Happel's theorem
(compare \cite{Buc}, Theorem 4.4.1 and \cite{Ha2}, Theorem 4.6, also
see \cite{Bel}, Theorem 6.9).

\begin{thm}
Let $R$ be a Gorenstein ring. Then the natural functor
$$\underline{F}: R\underline{\mbox{\rm -GProj}}\longrightarrow
D'_{\rm sg}(R)$$
is a triangle-equivalence.
\end{thm}

\noindent {\bf Proof.}\quad We have noted the following fact: set
$\mathcal{A}=R\mbox{-Mod}$, $\omega=R\mbox{-Proj}$, then
$_\omega\mathcal{X}=\mathcal{A}$ and
$\alpha(\omega)=R\mbox{-GProj}$. Hence by Theorem 2.1, we know that
to obtain the result, it suffices to show that
$\widehat{R\mbox{-GProj}}=R\mbox{-Mod}$. Assume ${\rm inj.dim}\;
_RR=d$. Then every projective module has injective dimension at most
$d$. Let $X$ be any $R$-module. Take an exact sequence
$$0\longrightarrow M \longrightarrow P^{d-1} \longrightarrow P^{d-2} \longrightarrow \cdots
\longrightarrow P^1 \longrightarrow P^0\longrightarrow X
\longrightarrow 0,$$ where each $P^i$ is projective. By
dimension-shift technique, we have, for each projective module $Q$,
${\rm Ext}_R^i(M, Q)\simeq {\rm Ext}_R^{i+d}(X, Q)=0$, $i\geq 1$.
Hence $M \in {^\perp R\mbox{-Proj}}$, and by Lemma 3.2, $M\in
R\mbox{-GProj}$. Hence, $X\in \widehat{R\mbox{-GProj}}$. Thus we are
done.\hfill $\blacksquare$

\subsection{}
In this subsection, we consider another self-orthogonal subcategory
$\omega'=R\mbox{-proj}$, the full subcategory of finite-generated
projective modules, of the category $\mathcal{A}=R\mbox{-Mod}$. From
the definition in {\bf 2.1}, it is not hard to see that
 \begin{align*} _{\omega'}\mathcal{X}=\{M \in
R\mbox{-Mod}\; &|\; \mbox{there exists an exact sequence } \\
\cdots \longrightarrow P^n \longrightarrow & P^{n-1}\longrightarrow
\cdots \longrightarrow P^1 \longrightarrow P^0\longrightarrow M
\longrightarrow 0, \; \mbox{each } P^n \in
R\mbox{-proj}\},\end{align*}
 and
\begin{align*}\mathcal{X}_{\omega'}=\{M \in R\mbox{-Mod}\; | & \; \mbox{there
exists
an exact sequence } \\
0 \longrightarrow M \longrightarrow  P^0
\stackrel{d^0}\longrightarrow P^1 \longrightarrow & \cdots
\longrightarrow P^n \stackrel{d^n}\longrightarrow
P^{n+1}\longrightarrow \cdots, \; \mbox{each } P^n\in R\mbox{-proj},
\; {\rm Coker}d^n \in {^\perp R\mbox{-proj}} \}.
\end{align*}
Set $\alpha(\omega')=R\mbox{-Gproj}$. Hence $R\mbox{-Gproj}$ is a
Frobenius category, whose projective-injective objects are exactly
contained in $R\mbox{-proj}$. Observe that $R\mbox{-Gproj}\subseteq
R\mbox{-GProj}$, and we have an induced inclusion of triangulated
categories $R\underline{\mbox{-Gproj}} \hookrightarrow
R\underline{\mbox{-GProj}}$.

Denote by $R\mbox{-mod}$ the full subcategory consisting of
finitely-presented modules. Let $R$ be a left-coherent ring. Observe
that in this case, $R\mbox{-mod}$ is an abelian subcategory of
$R\mbox{-Mod}$, and $R\mbox{-mod}={_{\omega'}\mathcal{X}}$ (compare
\cite{Au2}, p.41). Therefore, if $R$ is left-coherent, we have
\begin{align*}
R\mbox{-Gproj}=\{M \in R\mbox{-mod}\; | & \; \mbox{there exists
an exact sequence } \\
0 \longrightarrow M \longrightarrow  P^0
\stackrel{d^0}\longrightarrow P^1 \longrightarrow & \cdots
\longrightarrow P^n \stackrel{d^n}\longrightarrow
P^{n+1}\longrightarrow \cdots, \; \mbox{each } P^n\in R\mbox{-proj},
\; {\rm Coker}d^n \in {^\perp R\mbox{-proj}} \}.
\end{align*}

The following  observation is interesting.

\begin{lem} Let $R$ be a left-coherent ring. Then we have $R\mbox{\rm -GProj}\cap
R\mbox{\rm -mod}=R\mbox{\rm -Gproj}$.
\end{lem}

\noindent{\bf Proof.}\quad Let $M\in R\mbox{\rm -Gproj}$. Then we
have an exact sequence $ 0 \longrightarrow M \longrightarrow  P^0
\stackrel{d^0}\longrightarrow P^1 \longrightarrow  \cdots
\longrightarrow P^n \stackrel{d^n}\longrightarrow
P^{n+1}\longrightarrow \cdots $, where $P^i\in R\mbox{-proj}$ and
each ${\rm Coker} d^i\in {^\perp R\mbox{-proj}}$. Since each module
${\rm Coker} d^i$ is finitely-generated, and thus ${\rm Coker}
d^i\in {^\perp}R\mbox{-proj}$ implies that ${\rm Coker}d^i\in
{^\perp}R\mbox{-Proj}$ immediately. Thus we have $M\in
R\mbox{-GProj}$. Hence $R\mbox{-Gproj}\subseteq R\mbox{-GProj}\cap
R\mbox{-mod}$.

Conversely, assume that  $M\in R\mbox{-GProj}\cap R\mbox{-mod}$.
Then there exists an exact sequence $0 \longrightarrow M
\stackrel{\varepsilon}\longrightarrow P \longrightarrow
X\longrightarrow 0$, where $P\in R\mbox{-Proj}$ and $X\in
R\mbox{-GProj}$. By adding proper projective modules to $P$ and $X$,
we may assume that $P$ is free. Since $M$ is finitely-generated, we
may decompose $P=P^0\oplus P'^0$ such that $P^0$ is
finitely-generated and ${\rm Im}\varepsilon\subseteq P^0$. Consider
the exact sequence $0\longrightarrow M
\stackrel{\varepsilon}\longrightarrow P^0\longrightarrow
M'\longrightarrow 0$. We have $M'\oplus P'^0\simeq X$, and note that
$R\mbox{-GProj}\subseteq R\mbox{-Mod}$ is closed under taking direct
summands (by Proposition 5.1 in \cite{AR}, or \cite{EJ}), we deduce
that $M'\in R\mbox{-GProj}$. Observe that $M'\in R\mbox{-mod}$, and
we have $M'\in R\mbox{-GProj}\cap R\mbox{-mod}$. Applying the same
argument to $M'$, we can find an exact sequence $0\longrightarrow M'
\longrightarrow P^1 \longrightarrow M''\longrightarrow 0$ such that
$P^1$ is finitely-generated projective, and $M''\in
R\mbox{-GProj}\cap R\mbox{-mod}$. Continue this process, we can
derive a long exact sequence $0 \longrightarrow M \longrightarrow
P^0 \longrightarrow P^1\longrightarrow \cdots$. This is the required
sequence proving $M\in R\mbox{-Gproj}$.\hfill $\blacksquare$

Let $R$ be left-coherent. Set $\mathcal{A}'=R\mbox{-mod}$. The
relative singularity of $\mathcal{A}'$ with respect to $\omega'$ is
the usual singularity category of the ring $R$ (\cite{O1})
$$D_{\rm sg}(R)=D^b(R\mbox{-mod})/K^b(R\mbox{-proj}).$$

The following is read directly from Theorem 2.1.

\begin{prop}
Let $R$ be a left-coherent ring. The natural functor $\underline{F}:
R\underline{\mbox{\rm -Gproj}}\longrightarrow D_{\rm sg}(R)$ is a
fully-faithful triangle-functor.
\end{prop}

\begin{rem}
Consider the natural embedding $D^b(R\mbox{-mod})\hookrightarrow
D^b(R\mbox{-Mod})$, and observe that $K^b(R\mbox{-Proj})\cap
D^b(R\mbox{-mod})=K^b(R\mbox{-proj})$, and for any $P^\bullet\in
K^b(R\mbox{-Proj})$, $X^\bullet \in D^b(R\mbox{-mod})$, then any
morphism (inside $D^b(R\mbox{-Mod})$) from $P^\bullet$ to
$X^\bullet$ factors through an object of $K^b(R\mbox{-proj})$ (just
take a projective resolution $Q^\bullet\in K^{-, b}(R\mbox{-proj})$
of $X^\bullet$, then the brutally truncated complex $\sigma^{\geq
-n}Q^\bullet$, for large  $n$, is the required object). Now, It
follows that the natural induced functor $D_{\rm
sg}(R)\longrightarrow D'_{\rm sg}(R)$ is a full embedding (by
\cite{V1}, 4-2 Theorem). Finally, we have a commutative diagram of
fully-faithful triangle-functors
\[\xymatrix{
R\underline{\mbox{\rm Gproj}} \ar@{^{(}->}[r]^-{\underline{F}}
\ar@{^{(}->}[d] & D_{\rm sg}(R) \ar@{^{(}->}[d]\\
R\underline{\mbox{\rm -GProj}} \ar@{^{(}->}[r]^-{\underline{F}} &
D'_{\rm sg}(R). }\]
\end{rem}

\vskip 10pt

 A sufficient condition that the functor $\underline{F}$
in Proposition 3.5 is an equivalence is also that the ring $R$ is
Gorenstein. We need the following result.

\begin{lem}
Let $R$ be a Gorenstein ring. Then we have
$$R\mbox{\rm -Gproj}=\{M
\in R\mbox{\rm -mod}\; |\; {\rm Ext}_R^i(M, P)=0, \; P\in R\mbox{\rm
-proj}, \; i\geq 1\}.$$
\end{lem}

\noindent {\bf Proof.}\quad Just note that the left hand side is
equal to  $R\mbox{-mod}\cap {^\perp R\mbox{-Proj}}$. Then the result
follows from Lemma 3.2 and Lemma 3.4 directly. Let us remark that
the lemma can be also proved by the cotilting theory.\hfill
$\blacksquare$

Using Proposition 3.5 and Lemma 3.7 and applying a similar argument
as Theorem 3.3, we have the following result. Note that the result
was first shown by Buchweitz \cite{Buc} and its dual version was
shown independently by Happel in the finite-dimensional case
\cite{Ha2} (compare \cite{Bel}, Corollary 4.13 or \cite{CZ}, Theorem
2.5). A special case of the result was given by Rickard (\cite{Ric},
Theorem 2.1) which says that the singularity category of a
self-injective algebra is triangle-equivalent to its stable module
category (compare Keller-Vossieck \cite{KV}).

\begin{thm}
{\rm (Buchweitz-Happel)} \; Let $R$ be a Gorenstein ring. Then the
natural functor
$$\underline{F}: R\underline{\mbox{\rm -Gproj}}\longrightarrow
D_{\rm sg}(R)$$ is a triangle-equivalence.
\end{thm}

\subsection{} Let $T$ be a self-orthogonal object in any abelian category
$\mathcal{A}$. Set $\alpha(T)=\alpha({\rm add}\; T)$. We will relate
$\alpha(T)$ to the category of Gorenstein-projective modules over
the endomorphism ring.

\begin{thm}
Let $T$ be a self-orthogonal object, and let $R={\rm
End}_\mathcal{A}(T)^{\rm op}$. Then the functor ${\rm
Hom}_\mathcal{A}(T, -): \alpha(T) \longrightarrow R\mbox{\rm
-Gproj}$ is fully-faithful, and it induce a full embedding of
triangulated categories $\underline{\alpha(T)}\longrightarrow
R\underline{\mbox{\rm -Gproj}}$.
\end{thm}

Part of the theorem follows from  an observation of Xi (\cite{Xi},
Proposition 5.1), which we will recall. Let $T\in \mathcal{A}$ be
any object, $R={\rm End}_R(T)^{\rm op}$. Then we have  the functor
$${\rm Hom}_\mathcal{A}(T, -):\; \mathcal{A}\longrightarrow
R\mbox{-Mod}.$$ In general, it is not fully-faithful. However, it is
well-known that it induces an equivalence
$$ {\rm add}\; T\simeq
R\mbox{-proj},$$ in particular, the restriction  of ${\rm
Hom}_\mathcal{A}(T, -)$ to ${\rm add}\; T$ is fully-faithful.
Actually, we can define a larger subcategory, on which ${\rm
Hom}_\mathcal{A}(T, -)$ is fully-faithful. For this, recall that a
morphism $g:\; T_0\longrightarrow M$ with $T_0\in {\rm add}\;T$ is a
\emph{$T$-precover (= right $T$-approximation)} of $M$, if  any
morphism from $T$ to $M$ factors through $g$. Consider the following
full subcategory
\begin{align*}
{\rm App}(T):= \{M \in \mathcal{A}\; |\; &\mbox{there exists an
exact sequence } T_1\stackrel{f_1}\longrightarrow
T_0\stackrel{f_0}\longrightarrow M
\longrightarrow 0,\\
 &\; T_i\in{\rm add}\; T,\; f_0 \mbox{ is a }
T\mbox{-precovers},\; f_1: T_1\longrightarrow {\rm Ker}f_0 \mbox{ is
a } T\mbox{-precover}\}.
\end{align*}
For $M\in {\rm App}(T)$, such a sequence
$T_1\stackrel{f_1}\longrightarrow T_0\stackrel{f_0}\longrightarrow M
\longrightarrow 0$ will be called a \emph{$T$-presentation} of $M$.

The following result is contained in \cite{Xi} in slightly different
form.

\begin{lem}
The functor ${\rm Hom}_\mathcal{A}(T, -)$ induces a full embedding
of ${\rm App}(T)$ into $R\mbox{\rm -mod}$.
\end{lem}

\noindent {\bf Proof.}\quad The proof resembles the argument in
\cite{ARS}, p.102. Let $M\in {\rm App}(T)$ with $T$-presentation
$T_1\stackrel{f_1}\longrightarrow T_0\stackrel{f_0}\longrightarrow M
\longrightarrow 0$. Since $f_0$ and $f_1:T_1\longrightarrow {\rm
Ker}f_0$ are $T$-precovers, we have the following exact sequence of
$R$-modules
$${\rm Hom}_\mathcal{A}(T, T_1) \stackrel{{\rm Hom}_\mathcal{A}(T, f_1)}\longrightarrow {\rm Hom}_\mathcal{A}(T, T_0)
\stackrel{{\rm Hom}_\mathcal{A}(T, f_0)}\longrightarrow {\rm
Hom}_\mathcal{A}(T, M) \longrightarrow 0.$$ Recall the equivalence
${\rm Hom}_\mathcal{A}(T, -): {\rm add}\; T \simeq R\mbox{-proj}$.
Thus the left-hand side two terms in the sequence above are
finite-generated $R$-modules, and we infer that ${\rm
Hom}_\mathcal{A}(T, M)$ is a finite-presented $R$-module. Let $M'\in
{\rm App}(T)$ with $T$-presentation
$T'_1\stackrel{f'_1}\longrightarrow
T'_0\stackrel{f'_0}\longrightarrow M' \longrightarrow 0$. Given any
homomorphism of $R$-modules $\theta: {\rm Hom}_\mathcal{A}(T,M)
\longrightarrow {\rm Hom}_\mathcal{A}(T, M')$. Thus by a similar
argument as the comparison theorem in homological algebra, we have
the following diagram in $R\mbox{-mod}$
\[\xymatrix @C=55pt {
{\rm Hom}_\mathcal{A}(T, T_1) \ar[r]^-{{\rm Hom}_\mathcal{A}(T,
f_1)} \ar[d]^-{\theta_1} & {\rm Hom}_\mathcal{A}(T, T_0)
\ar[r]^-{{\rm Hom}_\mathcal{A}(T, f_0)} \ar[d]^-{\theta_0} & {\rm
Hom}_\mathcal{A}(T, M) \ar[r] \ar[d]^-\theta & 0 \\
{\rm Hom}_\mathcal{A}(T, T'_1) \ar[r]^-{{\rm Hom}_\mathcal{A}(T,
f'_1)} & {\rm Hom}_\mathcal{A}(T, T'_0) \ar[r]^-{{\rm
Hom}_\mathcal{A}(T, f_0)} & {\rm Hom}_\mathcal{A}(T, M') \ar[r] & 0.
}\]

Using the equivalence ${\rm add}\; T\simeq R\mbox{-proj}$ again, we
have morphisms $g_i: T_i\longrightarrow T_i'$  such that ${\rm
Hom}_\mathcal{A}(T, g_i)=\theta_i$, $i=0,1$. Thus $g_0\circ
f_1=f'_1\circ g_1$. Thus we have a unique morphism $g: M
\longrightarrow M'$ making the diagram commute

\[\xymatrix{
T_1 \ar[r]^-{f_1} \ar[d]^-{g_1} & T_0 \ar[r]^-{f_0} \ar[d]^-{g_0} &
M
\ar[r] \ar[d]^-g  & 0\\
T'_1 \ar[r]^-{f'_1} & T'_0 \ar[r]^-{f'_0} & M' \ar[r] & 0 .}\]

Now it is not hard to see that ${\rm Hom}_\mathcal{A}(T, g)=\theta$,
i.e., ${\rm Hom}_\mathcal{A}(T, -): {\rm App}(T) \longrightarrow
R\mbox{-mod}$ is full. We will omit the proof of faithfulness, which
is somehow the inverse of the above proof. \hfill $\blacksquare$

\noindent {\bf Proof of Theorem 3.9:} \quad Set $\omega={\rm add}\;
T$. First note that any epimorphism $f:T_0\longrightarrow M$ with
$T_0\in {\rm add}\; T$ and ${\rm Ker}f\in {T^\perp}$, is a
$T$-precover. This can be seen from the long exact sequence obtained
by applying ${\rm Hom}_\mathcal{A}(T, -)$ to the sequence $0
\longrightarrow {\rm Ker}f \longrightarrow T_0
\stackrel{f}\longrightarrow M\longrightarrow 0$. Thus we infer that
$_\omega \mathcal{X}\subseteq {\rm App}(T)$, and then
$\alpha(T)\subseteq {\rm App}(T)$. Thus ${\rm Hom}_\mathcal{A}(T,
-)$ is fully-faithful on $\alpha(T)$. What is left to show is that
for each $M\in \alpha(T)$, ${\rm Hom}_\mathcal{A}(T, M)\in
R\mbox{-Gproj}$. Take a complete $T$-resolution $T^\bullet=(T^n,
d^n)$ for $M$. Then the complex $P^\bullet={\rm Hom}_\mathcal{A}(T,
T^\bullet)$ is exact with its $0$-cocycle  ${\rm Hom}_\mathcal{A}(T,
M)$. Note that $P^\bullet$ is a complex of finitely-generated
projective $R$-modules. Note that we have an isomorphism of Hom
complexes ${\rm Hom}_\mathcal{A}(T^\bullet, T)\simeq {\rm
Hom}_R(P^\bullet, R)$, using the equivalence ${\rm
Hom}_\mathcal{A}(T, -): {\rm add}\; T\simeq R\mbox{-proj}$ and
noting that ${\rm Hom}_\mathcal{A}(T, T)=R$. However ${\rm
Hom}_\mathcal{A}(T^\bullet, T)$ is exact, hence we infer that
$P^\bullet$ is a complete resolution for ${\rm Hom}_\mathcal{A}(T,
M)$. Thus ${\rm Hom}_\mathcal{A}(T, M)\in R\mbox{-Gproj}$.

Note that ${\rm Hom}_\mathcal{A} (T, -)$ preserves short exact
sequences in $\alpha(T)$, and thus the composite $\alpha(T)
\longrightarrow R\mbox{-Gproj} \longrightarrow
R\mbox{\underline{-Gproj}}$ is a $\partial$-functor, which sends
${\rm add}\; T$ to zero. By Lemma 2.5, the induced functor
$\underline{\alpha(T)}\longrightarrow R\mbox{\underline{-Gproj}}$ is
a triangle-functor, the fully-faithfulness of which follows directly
from the one of ${\rm Hom}_\mathcal{A}(T, -):
\alpha(T)\longrightarrow R\mbox{-Gproj}$.\hfill $\blacksquare$

\section{Compact generators for Gorenstein-projective modules}

\subsection{} Let us begin with some notions. Let $\mathcal{C}$ be a
triangulated category with arbitrary (small) coproducts.  An object
$C\in \mathcal{C}$ is said to be \emph{compact}, if the functor
${\rm Hom}_\mathcal{C}(C, -)$ commutes with coproducts. Denote by
$\mathcal{C}^c$ the full subcategory  of $\mathcal{C}$ consisting of
compact objects, which is easily seen to be a thick triangulated
subcategory. The triangulated category $\mathcal{C}$ is said to be
\emph{compactly generated}, if there is a set $S$ of compact objects
such that there is no proper triangulated category containing $S$
and closed under coproducts \cite{Ne}.

Let $R$ be a ring. It is easy to see that the triangulated category
$R\mbox{\underline{-GProj}}$ has arbitrary coproducts, and the
natural embedding $R\mbox{\underline{-Gproj}} \longrightarrow
R\mbox{\underline{-GProj}}$ gives that
$R\mbox{\underline{-Gproj}}\subseteq
(R\mbox{\underline{-Gproj}})^c$.

We have our main result. Note that similar results were obtained by
Beligiannis (\cite{Bel2}, Theorem 6.7 and \cite{Bel3}, Theorem 6.6),
and Iyengar-Krause (\cite{IK}, Theorem 5.4 (2)) using different
methods in different setups. We would like to thank Beligiannis to
remark that one might find another proof of the following result
using Gorenstein-injective modules, and a suitable combination of
results and arguments in \cite{Bel3} and \cite{BR}.

\begin{thm}
Let $R$ be a Gorenstein ring. Then the triangulated category
$R\mbox{\underline{\rm -GProj}}$ is compactly generated, and its
subcategory of compact objects $(R\mbox{\underline{\rm -GProj}})^c$
is the additive closure of $R\mbox{\underline{\rm -Gproj}}$.
\end{thm}

Before proving Theorem 4.1, we need to recall some well-known facts
on the homotopy category of projective modules. Denote by $K_{\rm
proj}(R)$ the smallest triangulated category of $K(R\mbox{-Proj})$
containing $R$ and closed under coproducts. Denote by $K^{\rm
ex}(R\mbox{-Proj})$ the full subcategory of $K(R\mbox{-Proj})$
consisting of exact complexes. For each complex $P^\bullet\in
K(R\mbox{-Proj})$, there is a unique triangle
\begin{align*}
{\bf p}(P^\bullet) \longrightarrow P^\bullet \longrightarrow {\bf
a}(P^\bullet) \longrightarrow {\bf p}(P^\bullet) [1]
\end{align*}
with ${\bf p}(P^\bullet)\in K_{\rm proj}(R)$ and ${\bf
a}(P^\bullet)\in K^{\rm ex}(R\mbox{-Proj})$. Thus we have an exact
functor ${\bf a}: K(R\mbox{-Proj}) \longrightarrow K^{\rm
ex}(R\mbox{-Proj})$. Moreover, we have an exact sequence of
triangulated categories
\begin{align*}
0 \longrightarrow K_{\rm proj}(R) \stackrel{\rm inc}\longrightarrow
K(R\mbox{-Proj}) \stackrel{\bf a}\longrightarrow K^{\rm
ex}(R\mbox{-Proj}) \longrightarrow 0,
\end{align*}
where ``inc" denotes the inclusion functor (for details, see
\cite{Ke2} and compare \cite{Kr}, Corollary 3.9).

The following result is essentially due to J{\o}rgensen \cite{J}.

\begin{lem}
Let $R$ be a Gorenstein ring. Then the homotopy category
$K(R\mbox{\rm -Proj})$ is compactly generated, and its subcategory
of compact object is $K^{+, b}(R\mbox{\rm -proj})$.
\end{lem}

 \noindent {\bf Proof.}\quad To see the lemma, we need the results
 of J{\o}rgensen: let $R$ be a ring, recall the duality
 $*={\rm Hom}_R(-, R):\; R\mbox{-proj} \longrightarrow R^{\rm op}
 \mbox{-proj}$, which can be extended to another duality
 $*: \; K^{-}(R\mbox{-proj})\longrightarrow K^{+}(R^{\rm
 op}\mbox{-proj})$. By \cite{J}, Theorem 2.4, if the ring $R$ is
 coherent and every flat  $R$-module has finite projective dimension,
 then the homotopy category $K(R\mbox{-Proj})$ is
 compactly-generated, and then by \cite{J}, Theorem 3.2 (and its
 proof), the subcategory of compact objects is $K(R\mbox{-Proj})^c=\{P^\bullet
 \in K^{+}(R\mbox{-proj})\; |\; \mbox{the complex } (P^\bullet)^* \in K^{-,b}(R^{\rm
 op}\mbox{-proj})\}$. Note the following two facts: (1) for a
 Gorenstein ring $R$, every flat module has finite projective
 dimension by \cite{EJ}, Chapter 9, section 1; (2) for a Gorenstein
 ring $R$, we have an induced duality
 $*:\; K^{-, b}(R\mbox{-proj})\longrightarrow K^{+, b}(R^{\rm
 op}\mbox{-proj})$, which is because that the regular module has
 finite injective dimension. Combining the above two facts and
 J{\o}rgensen's results, we have the lemma. \hfill $\blacksquare$

Next result is also known, compare \cite{Buc}, Theorem 4.4.1.

\begin{lem}
Let $R$ be a Gorenstein ring. The following composite functor
\begin{align*}
 K^{\rm ex}(R\mbox{\rm -Proj}) \stackrel{Z^0}\longrightarrow R\mbox{\rm \underline{-GProj}}\end{align*}
 is a triangle-equivalence, where $Z^0$ is the functor of taking the zeroth
cocycyles.
\end{lem}

\noindent {\bf Proof.}\quad Note that since $_RR$ has finite
injective dimension, we infer that, by the dimension-shift
technique, every complex $P^\bullet \in K^{\rm ex}(R\mbox{-Proj})$,
its cocycles $Z^i$ lie in ${^\perp R\mbox{-Proj}}$, and furthur
$Z^i$ are Gorenstein-projective. Hence the above functor is
well-defined. Note that the functor is induced by the corresponding
functor of taking the zeroth cocycles $Z^0: C^{\rm
ex}(R\mbox{-Proj}) \longrightarrow R\mbox{\rm \underline{-GProj}}$,
and note that $Z^0$ is an exact functor between two exact
categories, preserving projective-injective objects. Hence the
induced functor $Z^0$ is a triangle-functor by \cite{Ha1}, p.23. The
proof of fully-faithfulness and denseness of $Z^0$ is same as the
argument in \cite{CZ}, Appendix (compare \cite{Bel}, Theorem 3.11).
Or, we observe that each exact complex $P^\bullet \in
K(R\mbox{-Proj})$ is a complete resolution (= totally-acyclic
complex in \cite{Kr}, section 7), and the result follows directly
from the dual of \cite{Kr}, Proposition 7.2. \hfill $\blacksquare$

\vskip 10pt

\noindent {\bf Proof of Theorem 4.1:}\quad We will see that the
result follows from the following result of
Thomason-Trobaugh-Yao-Neeman \cite{Ne}: let $\mathcal{C}$ be a
compactly generated and $S$ a subset of compact objects,
$\mathcal{R}$ the smallest triangulated subcategory which contains
$S$ and closed under coproducts, then the quotient category
$\mathcal{C}/\mathcal{R}$ is compactly generated, and every compact
objects in $\mathcal{C}/\mathcal{R}$ is a direct summand of $\pi(C)$
for some compact object $C$ in $\mathcal{C}$, where $\pi:
\mathcal{C} \longrightarrow \mathcal{C}/{\mathcal{R}}$ is the
quotient functor. To apply the theorem in our situation, by Lemma
4.2 we may put $\mathcal{C}=K(R\mbox{-Proj})$, and $S=\{R\}$, and
then $\mathcal{R}=K_{\rm proj}(R)$. Via the functor ${\bf a}$ and
the functor $Z^0$ in Lemma 4.3, we identify the quotient category
$\mathcal{C}/\mathcal{R}$ with $R\mbox{\underline{-GProj}}$. Hence
the triangulated category $R\mbox{\underline{-GProj}}$ is compactly
generated, every object $G$ in $(R\mbox{\underline{\rm -GProj}})^c$
is a direct summand of the image of the compact object in
$K(R\mbox{-Proj})$, and thus by Lemma 4.2 again, there exists
$P^\bullet \in K^{+, b}(R\mbox{-proj})$ such that $G$ is a direct
summand of $Z^0({\bf a}(P^\bullet))$.

Assume that $P^\bullet=(P^n, d^n)$, and take a positive number $n_0$
such that $H^n(P^\bullet)=0$, $n\geq n_0$. Consider the natural
distinguished triangle
\begin{align*}
\sigma^{\geq n_0}P^\bullet \longrightarrow P^\bullet \longrightarrow
\sigma^{< n_0}P^\bullet \longrightarrow (\sigma^{\geq
n_0}P^\bullet)[1],
\end{align*}
where $\sigma$ is the brutal truncation. Since
$\sigma^{<n_0}P^\bullet \in K^b(R\mbox{-proj}) \subseteq K_{\rm
proj}(R)$, we get ${\bf a}(\sigma^{<n_0}P^\bullet)=0$. Thus by
applying the exact functor ${\bf a}$ to the above triangle, we have
${\bf a}(P^\bullet)\simeq {\bf a}(\sigma^{\geq n_0}P^\bullet)$.
Applying the dimension-shift technique to the following exact
sequence and noting that the injective dimension of $_RR$ is finite
\begin{align*}
0 \longrightarrow Z^{n_0}(P^\bullet )\longrightarrow P^{n_0}
\stackrel{d^{n_0}}\longrightarrow P^{n_0+1} \longrightarrow \cdots
\longrightarrow P^n \stackrel{d^n} \longrightarrow P^{n+1}
\longrightarrow \cdots,
\end{align*}
we infer that $Z^{n_0}(P^\bullet)$ lies in ${^\perp}R\mbox{-proj}$,
and by Lemma 3.7, we have $Z^{n_0}(P^\bullet)\in R\mbox{-Gproj}$,
and thus it is not hard to see that ${\bf a}(\sigma^{\geq
n_0}P^\bullet)$ is a shifted complete resolution of
$Z^{n_0}(P^\bullet)$ (and in this case, ${\bf p}(\sigma^{\geq
n_0}P^\bullet)$ is the truncated projective resolution of
$Z^{n_0}(P^\bullet)$). Therefore $Z^0({\bf a}(\sigma^{\geq
n_0}P^\bullet))$ is the $n_0$-th syzygy of $Z^{n_0}(P^\bullet)$, and
thus it lies in $R\mbox{-Gproj}$. Hence $G$ is a direct summand of a
module in $R\mbox{-Gproj}$. This completes the proof. \hfill
$\blacksquare$

\vskip 10pt

\noindent {\bf Acknowledgement}\quad Some results in this paper
appeared in the second chapter of my Ph.D thesis, and I would like
to thank Prof. Pu Zhang for his supervision. I also would like to
thank Prof. Apostolos Beligiannis and Prof. Henning Krause very much
for their helpful remarks.

\bibliography{}

\end{document}